\input amstex
\documentstyle{amsppt}
\vsize=7.25in
\voffset=.25in
\hoffset=.25in
\magnification\magstep1
\nologo

\def\ZZ{{\Bbb Z}}

\def\RR{{\Bbb R}}

\def\fz{{F\times\Bbb{Z}}}
\def\boxx{\unskip \nopagebreak \hfill $\square$}

\topmatter
\title
Quasi-isometrically embedded subgroups\\
of Thompson's group $F$
\endtitle
\author
Jos\'e Burillo
\endauthor
\abstract
The goal of this paper is to construct quasi-isometrically embedded 
subgroups of Thompson's group $F$ which are isomorphic to 
$\fz^n$ for all $n$. A result estimating the norm of an element of Thompson's 
group is found. As a corollary, Thompson's group is seen to be an example 
of a finitely presented group which has an infinite-dimensional asymptotic 
cone.
\endabstract
\address
Dept. of Mathematics, Tufts University, Medford, MA 02155, U.S.A.
\endaddress
\email
jburillo\@emerald.tufts.edu
\endemail
\endtopmatter

\document
The interesting properties of Thompson's group $F$ have made it a
favorite object of study among group theorists and topologists. It
was first used by McKenzie and Thompson to construct finitely
presented groups with unsolvable word problems (\cite{5}). It is 
of interest also in homotopy theory in work related to homotopy 
idempotents, due to its universal conjugacy idempotent map $\phi$,
also used to see that $F$ is an infinitely iterated HNN extension. 
In \cite{1} Brown and Geoghegan found $F$ to be the
first torsion-free infinite-dimensional $FP_\infty$ group. Also,
$F$ contains an abelian free group of infinite rank, but it does not
admit a free non-abelian subgroup.

Many questions about $F$ are still open, in particular it is not known
whether $F$ is automatic, or what is its Dehn function ---although some
estimates have been found by Guba, who proves it is polynomial
in \cite{4}---. The amenability of $F$ is also unknown, fact that is of
considerable interest since both the affirmative or negative answer 
would provide counterexamples to open questions (see \cite{2}).

Questions about the geometric properties of $F$ have also been proposed. 
Bridson raised the question of whether $F$ and $\fz$ could be
quasi-isometric. In this paper we provide a partial answer to this question,
proving that $F$ admits a quasi-isometrically embedded subgroup isomorphic
to $\fz$. As a corollary, $F$ is the first example of a finitely
presented group whose asymptotic cone is infinite-dimensional (see \cite{3}).

There are in the literature several interpretations of $F$ that are
useful to study it. Cannon, Floyd and Parry provide two of these
interpretations, one as a group of homeomorphisms of the interval
$[0,1]$, and another as a group of isomorphisms of rooted binary trees.
In \cite{1} Brown and Geoghegan prove that $F$ is
isomorphic to certain group of piecewise linear homeomorphisms
of $\RR$, fact that will be used extensively in this paper. This construction
allows us to translate group-theoretical questions to the setting of
these homeomorphisms of $\RR$.

The organization of this paper is as follows: after a brief summary of 
results about $F$ in section 1, an estimate of the norm of an element
of $F$ is found in section 2. The last section is dedicated to state
and prove the results about the different subgroups of $F$.

The author would like to thank Z. Nitecki and J. Taback for useful
comments in the development of this work.

\heading
1. Generalities on Thompson's group $F$
\endheading

Thompson's group $F$ is the group defined by the following infinite
presentation:
$$
{\Cal P}=\left< x_k, k\ge 0 \,|\, x_i^{-1}x_jx_i=x_{j+1}, 
\text{ if }i<j \right>.
$$
Even though this is an infinite presentation, it follows
from the relations that the generators $x_k$, for $k>1$, are consequence
of $x_0$ and $x_1$. In fact, $F$ admits a finite presentation given by
$$
{\Cal F}=\left< x_0,\,x_1\,|\,[x_0x_1^{-1},x_0^{-1}x_1x_0],
[x_0x_1^{-1},x_0^{-2}x_1x_0^2]\right>
$$
(see \cite{1}). Throughout this paper, every time we refer to the word metric,
or the norm, or the distance in $F$, it will make reference to the finite
presentation $\Cal F$ of $F$.

It is a consequence of the presentation $\Cal P$ that
the map $\phi:F\longrightarrow F$ defined by $\phi(x_i)=x_{i+1}$ 
is a conjugacy idempotent, i.e., satisfies
$\phi^2(x)=x_0^{-1}\phi(x)x_0$, for all $x\in F$. Also, $\phi$ is
injective, mapping $F$ to the copy of itself generated by $x_k$, for $k\ge1$.
Moreover, this map is shown in \cite{1} to be universal among conjugacy
idempotents.

It is also seen in \cite{1} that the elements of $F$ admit a
unique normal form in the generators of $\Cal P$. 
From the relations in $\Cal P$ it is easy to see that any element of $F$
admits an expression of the form
$$
x_{i_1}^{r_1}x_{i_2}^{r_2}\ldots x_{i_n}^{r_n}
x_{j_m}^{-s_m}\ldots x_{j_2}^{-s_2}x_{j_1}^{-s_1}
$$
such that 
\roster
\item $r_1,\ldots,r_n,s_1,\ldots,s_m>0$, 
\item $i_1<i_2<\ldots<i_n$,
\item $j_1<j_2<\ldots<j_m$, and 
\item $i_n\ne j_m$. 
\endroster
This expression is not unique if we do not
require an extra condition: if for some $i$ both $x_i$ and $x_i^{-1}$ 
appear, then either $x_{i+1}$ or $x_{i+1}^{-1}$ must appear as well. For
otherwise there would be a subproduct of the form $x_i\phi^2(x)x_i^{-1}$
that could be replaced by $\phi(x)$. From this construction it is clear
that, given a word in the generators $x_k$ of $\Cal P$, using the relations 
and the extra condition, we can obtain the unique normal form, and the length 
of the word does not increase in this process. In other words,
the unique normal form is the shortest of all the words that represent
a given element in the generators of $\Cal P$. This fact will be used
later in this paper. 

Brown and Geoghegan show in \cite{1} that G admits an isomorphism with the group
of certain piecewise linear homeomorphisms of $\RR$. Let 
$$
f_k:\RR\longrightarrow\RR
$$
be the map defined by
$$
f_k(t)=\cases
t&\text{ if }t\le k,\\
2t-k&\text{ if }k\le t \le k+1,\\
t+1&\text{ if }t\ge k+1.
\endcases
$$
The group $G$ generated by the maps $f_k$, for all integers $k\ge 0$, and with 
right action, is isomorphic to Thompson's group, with $x_k$ identified 
with $f_k$. The right action means that the composition of two maps is
written the opposite way: the element $x_ix_j$ of $F$ is associated
with the element $f_if_j$ of $G$, which represents the map $f_j\circ f_i$
on $\RR$. From now on we will identify the groups $F$ and $G$.
This construction will be extremely useful in section 2.

Using this construction it is easy to see more properties of $F$: the
subgroup of $F$ generated by the elements $x_{2k}x_{2k+1}^{-1}$, for
$k\ge 0$, is a free abelian subgroup of infinite rank. To see that
two of these elements commute, observe that the map $f_{2n}f_{2n+1}^{-1}$
is the identity except in the interval $[2n,2n+2]$. Also, due to this fact,
it is clear that $x_0x_1^{-1}$ commutes with $x_k$ for $k\ge 2$, so
the subgroup generated by $x_0x_1^{-1}$, $x_2$ and $x_3$ is isomorphic
to $\fz$. In section 3 it will be proved that this subgroup
is nondistorted in $F$.

For an exhaustive survey of the properties of Thompson's groups (not
only $F$) and their geometric interpretations, see \cite{2}.

\heading
2. The estimate of the norm
\endheading

The geometric interpretation given by the maps $f_k$ provides a
method to compute the word metric of $F$ (with respect to the finite 
presentation $\Cal F$). For instance, to compute the norm of an
element $x$ of $F$, we can study the corresponding map $f$ in $G$:
we know that $f$ can be obtained as a composition
of the maps $f_0$ and $f_1$ and their inverses, and we only 
need to estimate how many occurrences of $f_0$ and $f_1$ and 
their inverses we need to obtain $f$. This can
be studied from properties of the graph of $f$. 

Given a point $(a,b)$ of the graph of $f$, with $b=f(a)$, since $f$ is
piecewise linear, we denote by $f'_+(a)$ and $f'_-(a)$ the right and left 
derivatives of $f$ at $a$. If $f'_+(a)\ne f'_-(a)$, we say that the point
$(a,b)$ is a breaking point of the graph of $f$.

To completely understand the maps in $G$ we need to study how multiplying 
by a generator affects a map. Let $f\in G$ be one of these piecewise
linear homeomorphisms. Then, the map $ff_i=f_i\circ f$ has a graph that can be
easily related to the graph of $f$. Since the map $f_i$ has slope 2
only on those points with $y$-coordinate in the interval $[i,i+2]$, the
graph of $ff_i$ is obtained by stretching the portion of the graph 
that has $y$-coordinate in $[i,i+1]$ to the interval $[i,i+2]$, and
all the graph is moved one unit up in all the points with $y\ge i+1$.
A point $(a,b)$ on the graph of $f$ with $b\in(i,i+1)$ appears now
as the point $(a,i+2(b-i))$, and the derivatives satisfy
$$
(ff_i)'_+(a)=2f'_+(a)\qquad\text{and}\qquad (ff_i)'_-(a)=2f'_-(a).
$$
Similarly, the map $ff_i^{-1}$ shrinks the interval $[i,i+2]$
down to $[i,i+1]$, and the derivatives get divided by 2. We will use
this facts in the proof of the norm estimate below.

The following lemma is an example of how this maps can be used to obtain 
group-theoretical properties:

\proclaim{Lemma 1} Let $f\in G$, and let $(a,b)$ be a point of the graph
of $f$. Assume that one of the two derivatives $f'_+(a)$ and $f'_-(a)$
is different from 1. Then the norm of $f$ can be bounded by
$$
|f|_G\ge \max\{1, a-2, b-2\}.
$$
In particular this applies to any breaking point of the graph of $f$.
\endproclaim

\demo{Proof} Since $|f|_G=|f^{-1}|_G$, we only need to prove that
$|f|_G\ge b-2$. And clearly we can assume $b\ge 3$.

Observe that in $f_0$ and $f_1$ the
highest point with a derivative different from 1 is the point $(2,3)$
in $f_1$, and further compositions with either $f_0$ or $f_1$ can 
only increase the $y$-coordinate by 1, and double the slope only
at a point with $y$-coordinate in $[0,3]$. So to achieve a derivative
different from 1 in $(a,b)$ one can start with $f_1$ and compose it with
$b-3$ generators more, at least. So one needs to compose at least $b-2$
generators to obtain a graph that has a derivative different from 1 
in $(a,b)$.\boxx
\enddemo

The next result is the estimate of the norm of an element of $G$ in
terms of the unique normal form.

\proclaim{Proposition 2} Let $f\in G$ be an element with normal form
$$
f=f_{i_1}^{r_1}\ldots f_{i_n}^{r_n} f_{j_m}^{-s_m}\ldots f_{j_1}^{-s_1}.
$$
Let $D=r_1+\ldots+r_n+s_1+\ldots+s_m+i_n+j_m$. Then
$$
\frac D6-2 \le |f|_G \le 3D.
$$
\endproclaim

\demo{Proof} Since $|f|_G=|f^{-1}|_G$, we can assume that $i_n>j_m$.

One of the inequalities is easy: rewrite the normal form in terms of
$f_0$ and $f_1$ using $f_i=f_0^{-i+1}f_1f_0^{i-1}$ to obtain a word
representing $f$ with only $f_0$, $f_1$ and their inverses. It is 
easy to see that the length of this word is less than $3D$.

It is also not difficult to see that $|f|_G\ge r_1+\ldots+r_n+s_1+\ldots+s_m$:
if $|f|_G<r_1+\ldots+r_n+s_1+\ldots+s_m$, there exists a word on $f_0$ 
and $f_1$ (and their inverses) that has length less than the unique 
normal form, contradicting the fact that the normal form is the shortest word 
for $f$. One could take this word and construct from it a normal form
that would be shorter than the unique one.

The last step in the proof is to prove that 
$$
|f|_G\ge \frac{i_n}2-2.
$$
Assume that $r_1+\ldots+r_n+s_1+\ldots+s_m<i_n/2$. If not, the inequality
follows from the previous paragraph. Consider the graph of the element
$$
g=f_{i_1}^{r_1}\ldots f_{i_n}^{r_n},
$$
the positive part of the normal form. Each one of these generators performs
a stretching of the graph (see above), the last one being a stretching
of the interval $[i_n,i_n+1]$ into $[i_n,i_n+2]$. So in the graph of $g$
there is a point $P$ with coordinates $(x,i_n+1)$ such that the two 
derivatives at this point are equal to $2^N$ where $N\ge r_n$. We want to
follow the movement down of $P$ after composing with all the inverses that 
appear in the normal form. The desired conclusion is that at the end, in
the graph of $f$, the point that corresponds to $P$ is a breaking point, or
else the function $f$ still has derivatives $2^N$ at this point.

Observe the effect that composing with $f_{j_m}^{-1}$ has on $P$. If 
$j_m=i_n-1$, in $gf_{j_m}^{-1}$ the point corresponding to $P$ is now a
breaking point: its derivatives are $2^N$ and $2^{N-1}$. Further 
compositions with any $f_i^{-1}$ will keep this point a breaking point,
since $i_n>j_m>\ldots>j_1$. If $j_m<i_n-1$, then in $gf_{j_m}^{-1}$
the point $P$ has just seen its $y$-coordinate decreased by one and the 
derivatives are still both $2^N$.

In the graph of $gf_{j_m}^{-s_m}$, the situation of the point $P$ depends
on the value of $s_m$:
\roster
\item If $s_m<i_n-j_m-1$, $P$ is not a breaking point and the derivatives
are both $2^N$.
\item If $s_m=i_n-j_m-1$, the last composition by $f_{j_m}^{-1}$ has made
$P$ a breaking point. 
\item If $s_m>i_n-j_m-1$, $P$ is now a breaking point whose $y$-coordinate 
is not an integer, and in any case it will remain a breaking point throughout,
so it will be a breaking point of $f$.
\endroster

The key to this argument is to observe that since $i_n>j_m>\ldots>j_1$,
a composition by an $f_i^{-1}$ cannot decrease the right derivative at $P$
without decreasing the left derivative. One of these compositions
either divides the left derivative by 2 without touching the right derivative, 
or it divides both derivatives by 2. And it divides both derivatives by 2
only after $P$ has been already made a breaking point.
So the conclusion is that in $f$, either the derivatives at
this point are still $2^N$ (if the $s_1,\ldots,s_m$ are small enough) or else
it is a breaking point. In any case one of the two derivatives at this point 
is not 1. We need now to compute the $y$-coordinate of this point
to apply Lemma 1.

In $g$ the $y$-coordinate of $P$ was $i_n+1$. Every application of an 
$f_i^{-1}$ may decrease the $y$-coordinate at most by 1. So the
$y$-coordinate is at least
$$
i_n+1-s_1-\ldots-s_m,
$$
but from our assumption, $r_1+\ldots+r_n+s_1+\ldots+s_m\le i_n/2$,
we conclude that the $y$-coordinate of this point is at least $i_n/2$.
From Lemma 1 it follows now that $|f|_G\ge \dfrac{i_n}2-2$.

Combining all inequalities (including $i_n>j_m$) we have
$$
\align
|f|_G&\ge\max\left\{r_1+\ldots+r_n+s_1+\ldots+s_m,\frac{i_n}2-2,
\frac{j_m}2-2\right\}\\
&\ge\frac{r_1+\ldots+r_n+s_1+\ldots+s_m+\dfrac{i_n}2+\dfrac{j_m}2-4}3\\
&\ge\frac D6-2,
\endalign
$$
and the proof is complete.\boxx
\enddemo

\heading
3. Quasi-isometrically embedded subgroups
\endheading

Recall that a map 
$$
F:X\longrightarrow Y
$$
between metric spaces is called a quasi-isometric embedding if there
exist constants $K,C>0$ such that
$$
\frac{d(x,x')}K-C\le d(F(x),F(x'))\le Kd(x,x')+C,
$$
for all $x,x'\in X$. If $G$ is a finitely generated group, and $H$ is
a finitely generated subgroup, then the fact that the inclusion is
a quasi-isometric embedding is equivalent to say that the distortion 
function
$$
h(r)=\frac1r \max\left\{|x|_H\,|\,x\in H, |x|_G\le r\right\}
$$
is bounded. If a subgroup is
quasi-iso\-metri\-cally embedded, then its own word metric is equivalent
to the metric induced by the word metric of the ambient group, and
the distortion is bounded.

Our goal is to prove that several subgroups of $F$ are quasi-isometrically
embedded. The first one and from which all the other ones will be
deduced, is the subgroup isomorphic to $\fz$ generated by the elements
$x_0x_1^{-1}$, $x_2$ and $x_3$. 

\proclaim{Theorem 3} The map
$$
\Phi:\fz\longrightarrow F
$$
defined by
$$
\Phi(x,t^k)=(x_0x_1^{-1})^k\phi^2(x)
$$
where $x\in F$ and $t$ is the generator of $\ZZ$, is a quasi-isometric 
embedding. 
\endproclaim

\demo{Proof} We need to prove that there exist constants $K,C>0$ such that
$$
\frac1K|(x,t^k)|_\fz-C\le |\Phi(x,t^k)|_F\le K|(x,t^k)|_\fz+C.
$$
Assume $x\in F$ has normal form
$$
x_{i_1}^{r_1}\ldots x_{i_n}^{r_n}x_{j_m}^{-s_m}\ldots x_{j_1}^{-s_1},
$$
and let 
$$
D=r_1+\ldots+r_n+s_1+\ldots+s_m+i_n+j_m.
$$
Then
$$
|(x,t^k)|_\fz=|k|+|x|_F,
$$
and, by the norm estimate,
$$
|k|+\frac D6-2\le
|(x,t^k)|_\fz\le |k|+3D.
$$
To compute the normal form of the element $\Phi(x,t^k)\in F$, assume
that $k\ge 0$. Observe that when $k<0$, we can just compute the norm of
$\Phi(x,t^k)^{-1}=\phi^2(x)^{-1}(x_0x_1^{-1})^{-k}=
(x_0x_1^{-1})^{-k}\phi^2(x^{-1})$.
The normal form is:
$$
\align
(x_0x_1^{-1})^k\phi^2(x)&=(x_0x_1^{-1})^k
x_{i_1+2}^{r_1}\ldots x_{i_n+2}^{r_n}
x_{j_m+2}^{-s_m}\ldots x_{j_1+2}^{-s_1}\\
&=x_0^kx_{i_1+k+2}^{r_1}\ldots x_{i_n+k+2}^{r_n}
x_{j_m+k+2}^{-s_m}\ldots x_{j_1+k+2}^{-s_1}
x_k^{-1}x_{k-1}^{-1}\ldots x_2^{-1}x_1^{-1}.
\endalign
$$
Then, 
$$
\frac{D+4k+4}6-2\le|\Phi(x,t^k)|_F\le3(D+4k+4).
$$
Combining these inequalities with the ones before, we obtain:
$$
\frac{|(x,t^k)|_\fz}{18}-2\le|\Phi(x,t^k)|_F\le 18|(x,t^k)|_\fz+48.
$$
which gives the desired result.\boxx
\enddemo

Once we know that $\fz$ quasi-isometrically embeds in $F$, we can extend 
this result to other subgroups:

\proclaim{Corollary 4} Ths subgroup generated by the elements
$$
x_0x_1^{-1},x_2x_3^{-1},\ldots,x_{2n-2}x_{2n-1}^{-1},x_{2n},x_{2n+1}
$$
is isomorphic to $\fz^n$ and it is quasi-isometrically embedded.
\endproclaim

Also, the free abelian subgroup $\ZZ^n$ in $\fz^n$ is also 
quasi-isometrically embedded. This implies that the asymptotic cone
of $F$ contains copies of $\RR^n$ for arbitrary $n$. For definitions
and properties of asymptotic cones see \cite{3}.

\proclaim{Corollary 5} The asymptotic cone of $F$ is infinite-dimensional.
\endproclaim

\enddocument